\documentclass[12pt, openany, amssymb, psamsfonts]{amsart}
\usepackage{mathrsfs,comment}
\usepackage[usenames,dvipsnames]{color}
\usepackage[normalem]{ulem}
\usepackage{url}
\usepackage[all,arc,2cell]{xy}
\UseAllTwocells
\usepackage{enumerate}
\usepackage{Template}

\usepackage{setspace}
\title{Spectral orthogonality of some special flows}

\author{Mingcheng Sheng}

\date{Fall 2025}

\usepackage[textwidth=13.5cm, centering]{geometry}

\begin{document}

\begin{abstract}
    In this paper, we study the spectral orthogonality problem for special flows built over irrational rotations under two different types of roof functions:
    
    1) the roof functions are real analytic.
    
    2) the roof functions are piecewise $C^1$ with one discontinuity. These flows are also known as von-Neumann flows.

    We show that if $\{T^f_\alpha\}$ is as in 1) and weak mixing, then for a $G_\delta$ dense set of $\beta$, we have that $\{T^f_\beta\}$ is weak mixing and is spectrally orthogonal to $\{T^f_\alpha\}$. On the other hand, if $\{T^f_\alpha\}$ is as in 2), then for a full measure set of $\beta$, the flows $\{T^f_\alpha\}$ and $\{T^f_\beta\}$ are spectrally orthogonal.
\end{abstract}
\maketitle

\section{Introduction}
Given two continuous probability preserving dynamical flows $\mathcal{T}=(T_t)_{t\in\mathbb{R}}$ on $(X,\mu)$ and $\mathcal{R}=(R_t)_{t\in\mathbb{R}}$ on $(Y,\nu)$, there are several ways to measure how different these systems are. The most classical is the {\em isomorphism} relation, where the two systems are isomorphic if there exists an invertible measure preserving transformation $S : X \to Y$ such
that $S \circ T_t = R_t \circ S$ for every $t \in \mathbb{R}$. A stronger notion is {\em disjointness} introduced by Furstenberg \cite{Furstenberg1967}. Two systems are disjoint (in the sense of Furstenberg) if the only {\em joining} between them is the product measure. An even stronger notion is {\em spectral orthogonality} or spectral disjointness. Two systems are spectrally disjoint if their spectral measures (see Section \ref{Spectral Orthogonality}) are orthogonal. It follows that spectral orthogonality implies (Furstenberg) disjointness which implies non-isomoprhism. 
This paper is concerned with studying the spectral orthogonality for some classes of special flows on the two torus that come from reparametrizations of linear toral flows. Recall that for the classical linear flows on $\mathbb{T}^2$ the isomorphism problem is well understood. It then follows that two such flows are spectrally orthogonal if and only if their rotation numbers are linearly independent over $\mathbb{Q}$. Note that this condition is also equivalent to Furstenberg disjointness. A more complicated case arises from smooth reparametrizations of linear flows. These flows were first studied by Kolmogorov \cite{Kol1953} who showed that if the rotation number of the linear flow is Diophantine (see \ref{Continued fractions}), then the reparametrized flow is smoothly conjugated to the corresponding linear flow and so has a pure point spectrum. However, the case of a Liouvillian rotation number is fundamentally different. In 1967, Shklover \cite{shklover1967} proved that there exist weak mixing analytic reparametrizations, which implies that the corresponding spectrum is continuous. The main goal of this paper is to study the spectral orthogonality for the class of weakly mixing reparametrizations of Liouvillean flows. Considering a Poincare section, these flows can be viewed as special flows over irrational rotations(see the introduction of \cite{FAYAD2002} for a short explanation). In the case where the roof functions have symmetric logarithmic singularities, Berk and Kanigowski proved in \cite{BerkKanigowski2021} that the rescaled special flows are spectrally orthogonal for a full measure set of rotation numbers. In contrast, we studied the case where the roof functions are analytic. Denote by $T_\alpha^f$ the special flow over a rotation $R_\alpha(x)= x+\alpha$ under a roof function $f$(see \ref{Special flows over irrational rotation}).

Our main results are as follows. 
\begin{thm}\label{main thm}

    Let $T_\beta^g$ be a special flow with a smooth roof function $g$, and
let $f$ be a real analytic function that is not a trigonometric polynomial. Then there exists a dense
$G_\delta$ set of irrational rotation number $\alpha$ such that the special flow $T^f_\alpha$ is spectrally disjoint from $T_\beta^g$.
\end{thm}
The requirement that $f$ is analytic but not a trigonometric polynomial ensures that we can construct weak-mixing flows. Otherwise, the following well known theorem shows that the flow can not be weak-mixing.
\begin{thm} \cite{Kol1953}
    $T_\beta^f$ is weak mixing implies that $f$ is not a trigonometric polynomial.
\end{thm}

The orthogonality problem has been studied for several classes of flows and maps. For example, Chaika\cite{Chaika2012} showed in 2012 that every ergodic transformation is
disjoint(in the sense of Furstenberg) from almost every interval 
exchange transformation. Also, Krzysztof Fraczek, Adam Kanigowski, and Corinna Ulcigrai\cite{kanigowski2025} showed that for almost every pair of IETs (Interval Exchange Transformation), the corresponding special flows are spectrally orthogonal. More results on spectral orthogonality and Furstenberg disjointness can be found in \cite{KanigowskiLemanczykUlcigrai2020},\cite{KanigowskiLemanczykRadziwill2021},\cite{Lemanczyk2022}, etc.

Our method also applies to another type of special flows called the Von Neumann flows, see chapter \ref{Von flow Chapter} for the definition. In this case, we show:
\begin{thm}\label{Von main thm}
    For every $\alpha\in[0,1]\setminus\mathbb{Q}$, there exists a full measure set $S_\alpha$, $m(S_\alpha)=1$, such that for every $\beta\in S_\alpha$,  the flows $\mathcal{T} = T^f_\alpha$ and $\mathcal{R} = T^f_\beta$ are spectrally orthogonal.
\end{thm}

The isomorphism problem of Von Neumann flows was first studied by Kanigowski and Solomko \cite{KANIGOWSKI_SOLOMKO_2018}. They showed that every two Von Neumann flows with the same rotation number but under roof functions of different slope are always non-isomorphic. In contrast, we establish a much stronger conclusion of spectral orthogonality in a more subtle setting where two roof functions share the same slope. The trade-off is that our result applies only to a full measure set of system parameters.
\subsection{Strategy and organization of the rest of the paper}

Our main tool to prove spectral orthogonality is a criterion from \cite{kanigowski2025}, which roughly says that if there exists a time sequence $\{t_n\}$ that serves as a rigidity sequence for one flow but a mixing sequence for the other, then the two flows are spectrally orthogonal. We will later establish that the sequence generated by the convergents of the base rotation number always constitutes a rigidity sequence. The problem is thus reduced to showing that a large set of special flows exhibits mixing behavior along the same sequence of convergents.

Let us briefly describe the content of the following
Sections. In section \ref{Basic Defn} we first introduce the background definitions, and then state our main criterion \ref{main criterion}. In section \ref{mixing subsequence chap} we reduce the problem of constructing the rigid and mixing subsequences to the estimation of Birkhoff sums of $f'$. Then in the next two sections \ref{sec basic estimates} and \ref{sec Estimation of the birkhoff sum}, we estimate the Birkhoff sum for the analytic roof functions. In section 6, we construct a $G_\delta$ set of $\alpha$ such that $T_\alpha^f$ mixes along a subsequence of the rigid sequence of $T_\beta^f$, thus completing the proof of Theorem \ref{main thm}. Finally, sections \ref{Von flow Chapter} and \ref{sec beta full measure} are devoted to the proof of Theorem \ref{Von main thm}. The proof of this result relies on the same machinery developed in sections \ref{Basic Defn} and \ref{mixing subsequence chap}.

\section{Basic Definitions}\label{Basic Defn}
\subsection{Special flows over irrational rotation}\label{Special flows over irrational rotation}
\textit{Birkhoff sums and special flows. }

Let $\mathbb{T}=\mathbb{R}/\mathbb{Z}$ be the one dimensional torus. For an ergodic automorphism $T$ on $\mathbb{T}$ and a strictly positive function $f\in L^1(\mathbb{T})$, we define a cocycle $\mathbb{Z}\times\mathbb{T}\to\mathbb{R}$ by: 
\begin{equation}
    S_n(f)(x) =\begin{cases}
      \sum_{0\le i<n} f(T^ix) &\text{ if } n\ge0\\
      0 &\text{ if } n=0\\
\sum_{n\le i<0} f(T^ix) &\text{ if } n<0
    \end{cases}
\end{equation}

The special flow $\mathcal{T}=\{T_t^f\}_{t\in\mathbb{R}}$ over the base automorphism $T : \mathbb{T} \to \mathbb{T}$ under the roof function $f$ is the unit vertical speed flow on the region $\mathbb{T}^f := \{(x,y) \in \mathbb{T} \times\mathbb{R} : 0 \le y < f(x)\}$, where we identify each point in the graph of the form $(x,f(x)),x\in \mathbb{T}$ with the base point
$(T(x),0)$. More precisely, for $x\in\mathbb{T}$, define $n:=N(x,t)\in\mathbb{Z}$ the unique integer number such that $S_n(f)(x) \le r+ t<S_{n+1}(f)(x)$. Then $(T^f_
t )_{t\in\mathbb{R}}$ acts on $\mathbb{T}^f$ so that
\begin{equation*}
    T^f_t (x,r) = (T^nx,r+ t-S_n(f)(x)).
\end{equation*}

In this paper, to simplify the
exposition of our proofs, we assume that $\int_{\mathbb{T}} f=\int_{\mathbb{T}} g=1$ and that 
\begin{equation*}
    \frac{1}{2}\le f\le\frac{3}{2}
\end{equation*}
which does not cause any loss of generality in the propositions we shall state. For general functions, we only need to change the constant accordingly.

\subsection{Continued fractions}\label{Continued fractions}
For an irrational $\alpha \in \mathbb{R}\setminus\mathbb{Q}$, $\alpha=[a_0;a_1, a_2, \cdots]$ is called the continued fraction
expansion of $\alpha$ and the rational numbers $(p_n/q_n)_{n\in\mathbb{N}}=[a_0;a_1\cdots,a_n]$ are called the convergents of
$\alpha$(see \cite{EinsiedlerWard2011}). The sequence $(q_n)_{n\in\mathbb{N}}$ is given by the recursive formulas
\begin{equation}\label{recursive formula continued fraction}
    q_0 = 1, q_1 = a_1, q_n = a_nq_{n-1} +q_{n-2}
\end{equation}
and the following two inequalities hold
\begin{equation}\label{4}
    \|q_{n-1}\alpha\|\le\|k\alpha\|, \quad \forall k<q_n
\end{equation}
where $\|x\|$ denotes the distance between $x$ and its closest integer, and 
\begin{equation}\label{5}
    \frac{1}{2q_{n+1}}\le \frac{1}{q_n+q_{n+1}}\le\|q_n\alpha\|<\frac{1}{q_{n+1}} \quad \forall n\in\mathbb{N}
\end{equation}
\begin{defn}(Liouville Number)
 An irrational number $x$ is a Liouville number if for every positive integer $n$, there exist positive integers $p$ and $q$ such that
    \begin{equation*}
        0<\bigg|x-\frac{p}{q}\bigg|<\frac{1}{q^n}
    \end{equation*}
\end{defn}
\begin{defn} (Diophantine number)
 An irrational number $x$ is Diophantine if it is not Liouville.
\end{defn}
\subsection{Spectral Orthogonality}\label{Spectral Orthogonality}
We now recall basic definitions from spectral theory. For a measure preserving dynamical flow $\mathcal{T}=\{T_t\}_{t\in\mathbb{R}}$ acting on the space $(X,\mathcal{B},\mu)$, the Koopman operators, denoted by $\{U^{\mathcal{T}}_{t}\}_{t\in \mathbb{R}} $, are a group of unitary operators in $L^2(X,B,\mu)$ given by $f \to f \circ T_t$. A subspace $\overline{\text{span}}\{f \circ T_t; t\in\mathbb{R}\}$ is a {\em cyclic space} of an element $f\in L^2(X,B,\mu)$. Due to
the Bochner-Herglotz theorem, there exists a finite measure $\sigma_f$ on $\mathbb{R}$ such that its Fourier transform satisfies
\begin{equation*}
    \hat{\sigma}_f(t)=\langle f\circ T_{-t},f\rangle \text{ for every } t\in\mathbb{R}
\end{equation*}
We say that $\sigma_f$ is the spectral measure of $f$. A maximal spectral type of $T$ is a spectral measure $\sigma_T$ such
that for every $f\in L^2(X,B,\mu)$, the spectral measure $\sigma_f$ is absolutely continuous with respect to $\sigma_T$. 
\begin{defn}(Spectral orthogonality)
    Two measure-preserving dynamical systems, $(X,\mu,T)$ and $(Y,\nu,S)$, are defined as spectrally orthogonal if the maximal spectral types of their associated Koopman operators, $U_T$ and $U_S$, are mutually singular.
\end{defn}
 Before stating the criterion of spectral orthogonality, let us first recall the definition of a rigidity sequence and a mixing sequence:
\begin{defn}(Rigidity sequence)
   An automorphism $T$ of a probability Borel space $(X,\mathcal{B},\mu)$ is called
rigid if there exists an increasing sequence $(h_n)_{n\in\mathbb{N}}$ of natural numbers such that
\begin{equation*}
    \lim_{n\to+\infty}\mu(A\triangle T^{h_n}A) = 0 \textit{ for every } A\in \mathcal{B}.
\end{equation*}
The sequence $(h_n)_{n\in\mathbb{N}}$ is called a rigidity sequence for $T$.
\end{defn}
\begin{defn}(Mixing sequence)
    A measure preserving flow $(T_t)_{t\in\mathbb{R}}$ in $(X,\mu)$ is mixing along an increasing sequence $(d_n)_{n\in \mathbb{N}}$ of real numbers if
    \begin{equation*}
        \lim_{n\to\infty}\langle U^{d_n}g_1,g_2\rangle= \langle g_1,1\rangle\langle 1,g_2\rangle \textit{ for every } g_1,g_2 \in L^2(X,\mu),
    \end{equation*}
\end{defn}
Now we state the main criterion:
\begin{thm}(Singularity Criterion via rigidity and exponential tails)\cite{kanigowski2025} Theorem 5. Chapter 4.2 \label{main criterion}

 Let $g : X \to\mathbb{R}_{>0}$ be an integrable roof function with $\inf_{x\in X} g(x) > 0$. Suppose that there exists a
rigidity sequence $(h_n)_{n\in\mathbb{N}}$ for $T$. Assume that $(T^g_t )_{t\in\mathbb{R}}$ on $(X^g,\mu^g)$ is a special flow such that there exist $(A_n)_{n\in \mathbb{N}}$ Borel sets with $\mu(A_n) \to1 $ as $n \to \infty$, and $(d_n)_{n\in\mathbb{N}}$ such that
$(S_{h_n}(g)(x)-d_n)_{n\in\mathbb{N}}$ has exponential tails, i.e. there exist two positive constants $C$ and $b$ such that
\begin{equation}
   \mu({x\in A_n : |S_{h_n}(g)(x)-d_n|\ge t}) \le Ce^{-bt} \textit{ for all } t\ge0 \textit{ and } n\in\mathbb{N}.
\end{equation} 
If $(U_t)_{t\in\mathbb{R}}$ is a measure-preserving flow on $(Y,\nu) $ that mixes along the sequence $(d_n)_{n\ge1}$,
then the flows $(T^g_t )_{t\in\mathbb{R}}$ and $(U_t)_{t\in\mathbb{R}}$ are spectrally orthogonal.
\end{thm} 

\section{Rigid and Mixing subsequence}\label{mixing subsequence chap}
The following classical consequence of the Denjoy-Koksma inequality shows that we can take the sequence $d_n$ in \ref{main criterion} to be $d_n=q_n\int g$.
\begin{lem}\label{lemma on rigid sq}\cite{Herman1979}
    If $g\in C^1$ and $\{p_n/q_n\}$ is the set of convergents of $\beta$, then $|S_{q_n}(x)-q_n\int g\ dm_L|_{C^0}\to 0$ uniformly. 
\end{lem}
Therefore, our attention can be restricted exclusively to the mixing subsequence. The main idea to prove mixing along a sequence is to estimate the behavior of the Birkhoff sum. We would like to show that for arbitrarily small cubes, the Birkhoff sum on the cube is \textit{stretching} evenly to infinity, and therefore the mass of the cube will spread evenly across the whole space. In this chapter, we will make this idea precise. 

\begin{defn}(Uniform stretch)\cite{FAYAD_2002} Let $\epsilon>0$ and $K>0$. We say that a real function $g$ on
an interval $[a,b]$ is $(\epsilon,K)$-uniformly stretching on $[a,b]$ if
\begin{equation*}
    \sup_{[a,b]}g- \inf_{[a,b]} g\ge K,
\end{equation*}
and if for any $u$ and $v$ such that
\begin{equation*}
  \inf_{[a,b]} g\le u< v\le \sup_{[a,b]}g  
\end{equation*}
the set
\begin{equation*}
    I_{u,v} = \{x\in [a,b]:u\le g(x)\le v\}
\end{equation*}
has Lebesgue measure
\begin{equation*}
    (1-\epsilon)\frac{v-u}{g(b)-g(a)}(b-a)\le \lambda(I_{u,v})\le (1+\epsilon)\frac{v-u}{g(b)-g(a)}(b-a)
\end{equation*}
We also write $g[a,b]$ for $\sup_{x\in[a,b]} |g(x)| - \inf_{x\in[a,b]} |g(x)|$. 
\end{defn}
For any $x \in I$ and $m > 0$, we
introduce the notation
\begin{equation*}
    \underline{N}(I,m):=\min_{x\in I} N(x,m), \quad \overline{N}(I,m):=\max_{x\in I} N(x,m)
\end{equation*}
for the maximum and minimum of the function $x\to N(x,t)$ on the interval $I$ respectively. ($N(x,t)$ is defined in section \ref{Special flows over irrational rotation})
\begin{prop}\label{criterion of mixing by uniform stretch}\cite{FAYAD_2002}
    Let $\{T_\alpha^f\}$ be a special flow over the irrational rotation $R_\alpha$ with a ceiling function $f>0$. Assume that for a sequence $t_n \to\infty$ we have that there
exists a sequence $(\mathcal{P}_n)_{n\ge1}$ of finite collections of disjoint intervals (called partial partitions), and positive functions $\epsilon(n),k(n)$ satisfying: 
\begin{equation*} 
    \text{Leb}(\mathcal{P}_n) \to 1; \quad\sup_{I\in\mathcal{P}_n}\text{Leb}(I)\to 0\quad \text{as} \quad n \to \infty
\end{equation*}
\begin{equation*} 
    \epsilon(n)\xrightarrow[n\to\infty]{}0; \quad k(n)\xrightarrow[n\to\infty]{}\infty
\end{equation*}
and for any $I\in\mathcal{P}_n$, $r \in [\underline{N}(I, t_n), \overline{N}(I, t_n)]$, 
\begin{equation*} 
    S_r(f)(x) \textit{ is $(\epsilon(n),k(n))$- uniformly stretching on } I.
\end{equation*}
Then $\{T_\alpha^f\}$ mixes along $t_n$.
\end{prop}
This is in fact a simplified version of the theorem in \cite{FAYAD_2002}, where the author deals with a two-dimensional rotation on $\mathbb{T}^2$.

\begin{prop}(Criterion of mixing along a sequence)\label{Criterion of mixing along a sequence}\cite{kanigowski2025} Lemma 6.5. 

Let $\{T_\alpha^f\}$ be a special flow over the irrational rotation $R_\alpha$ with a ceiling function $f>0$ of class at least $C^2$. Assume that for a sequence $t_n \to\infty$ we have that there
exists a sequence $(\mathcal{P}_n)_{n\ge1}$ of finite collections of disjoint intervals (called partial partitions) satisfying: 
\begin{equation} \tag{M} \label{eq:M}
    \text{Leb}(\mathcal{P}_n) \to 1 \quad \text{as} \quad n \to \infty;
\end{equation}
\begin{equation} \tag{S1} \label{eq:S1}
    \min_{I \in \mathcal{P}_n} \min_{r \in [\underline{N}(I, t_n), \overline{N}(I, t_n)]} \min_{x \in I} |S_r(f')(x)| |I| \to \infty \quad \text{as} \quad n \to \infty;
\end{equation}

\begin{equation} \tag{S2} \label{eq:S2}
    \max_{I \in \mathcal{P}_n} \frac{\max_{0 \le r \le 2t_n} \max_{x \in I} |S_r(f'')(x)| |I|}{\min_{r \in [\underline{N}(I, t_n), \overline{N}(I, t_n)]} \min_{x \in I} |S_r(f')(x)|} \to 0 \quad \text{as} \quad n \to \infty.
\end{equation}
then $\{T_\alpha^f\}$ mixes along $t_n$.
    
\end{prop}
\begin{rem}
    The original version of this lemma in \cite{kanigowski2025} has one additional condition: $T^jI\cap End(T) = \varnothing $. $End(T)$ denotes the set of discontinuities of the interval exchange transform. However, this condition is automatically satisfied in our case, because there is no discontinuity for the irrational rotation on $\mathbb{T}$.
\end{rem}

We will also use the following estimate on $N(x,t)$. (Lemma 1.2 \cite{FAYAD_2002}) 
\begin{lem}\label{N(x,t) estimation}
    If $f$ is analytic and satisfies the condition in section \ref{Special flows over irrational rotation}, then for $t$ large enough, and any $x\in \mathbb{T}$,
    \[N(x,t)\in [t/2,2t]\]
\end{lem} 
\section{Basic estimates}\label{sec basic estimates}
In this chapter, we collect some inequalities that are useful for estimation of the Birkhoff sum. Assume that the analytic function $f$ has the Fourier expansion $f(x)=\hat{f}_0+\sum_{k\neq 0}\hat{f}_ke^{i2\pi kx}$. Then we compute 
\[S_m(f)(x)=m\hat{f}_0+\sum_{k\neq 0}\hat{f}_ke^{i2\pi kx}\frac{1-e^{i2\pi km\alpha}}{1-e^{i2\pi k\alpha}}\]
We will adopt the notation in \cite{FAYAD_2002} and write

\begin{equation}
    X(m,k):=\frac{1-e^{i2\pi km\alpha}}{1-e^{i2\pi k\alpha}}=\sum_{j=0}^{m-1}e^{i2\pi jk\alpha}
\end{equation}

 Now we discuss strategies to bound $|X(m,k)|$ for different ranges of $m,k$. Since $X(m,-k)=\overline{X(m,k)}$, all bounds in this section assume that $k\in \mathbb{N}^*$. Let $q_n< m< q_{n+1}$.  
 Starting at the lower bounds, by definition, we have the trivial bound
    \begin{equation}\label{1}
        \forall k\in\mathbb{N^*}, \,m\in\mathbb{N}, \quad |X(m,k)|\le m
    \end{equation}
    
Another bound is 
\begin{equation*}
    |X(m,k)|\le\frac{2}{|1-e^{i2\pi k\alpha}|}.
\end{equation*}

When $0\le u\le 1/2$, we have the inequality $\sin{\pi u}\ge 2u$, so 
\begin{equation}
    \frac{2}{|1-e^{i2\pi k\alpha}|}=\frac{1}{\sin{\pi\|k\alpha\|}}\le\frac{1}{2\|k\alpha\|}
\end{equation}
    
    When $k$ is comparable with $q_n$, by \ref{4} and \ref{5} we have the following
    \begin{equation}\label{2}
        \forall n\in\mathbb{N},\, k<q_n,\, m\in\mathbb{N}, \quad |X(m,k)|\le q_n
    \end{equation}

When $k\in(jq_n,(j+1)q_n)$, where $j<q_{n+1}/4q_n$, we can write $k$ in the form $k=jq_n+q$, $0<q<q_n$, then 
\begin{equation*}
\|k\alpha\|=\|jq_n\alpha+q\alpha\|\ge\|q\alpha\|-\|jq_n\alpha\|\ge \|q_{n-1}\alpha\|-j\|q_n\alpha\|\ge\frac{1}{2q_n}-\frac{j}{q_{n+1}}\ge \frac{1}{4q_n}
\end{equation*}
Therefore, (remark, this only works for $k\ne jq_n$)
    \begin{equation}\label{3}
        \forall n\in\mathbb{N},\, k\in(jq_n,(j+1)q_n), \text{with } j<q_{n+1}/4q_n,\, m\in\mathbb{N}, \quad |X(m,k)|\le 2q_n
    \end{equation}
Finally, for any $m\le q_{n+1}/2$, we also have the following estimates from \cite{FAYAD_2002}:
   \begin{equation}\label{7}
        |\arg X(m,q_n)|\le \pi\frac{m}{q_{n+1}}
    \end{equation}
    \begin{equation}\label{6}
        |X(m,q_n)|\ge \frac{2}{\pi}m
    \end{equation}
\enlargethispage{3.5\baselineskip}    
\section{Estimation of the Birkhoff sum}\label{sec Estimation of the birkhoff sum}

To apply criterion \ref{Criterion of mixing along a sequence}, we need to estimate the lower bound of $S_m'$ and the upper bound of $S_m''$. We first estimate  $S^\prime_m$, assuming $q_n<m<q_{n+1}/2$. By direct computation and using the fact that $X(m,-k)=\overline{X(m,k)}$, we can split $S_mf'$ into the main term and other error terms:
\begin{equation*}
    \begin{split}S_m(f^\prime)(x)&=\sum_{k\neq 0}2\pi ik\cdot \hat{f}_ke^{i2\pi kx}X(m,k)\\&= \underbrace{2\pi iq_n \cdot (\hat{f}_{q_n}e^{i2\pi q_nx}-\hat{f}_{-q_n}e^{-i2\pi q_nx})|X(m,q_n)|}_{S_1(\text{main term})}\\&+\underbrace{
    \begin{aligned}
        2\pi iq_n \bigg[ &\hat{f}_{q_n}e^{i2\pi q_nx}\bigg(X(m,q_n)-|X(m,q_n)|\bigg)\\&- \hat{f}_{-q_n}e^{-i2\pi q_nx}\bigg(X(m,-q_n)-|X(m,-q_n)|\bigg)\bigg]
    \end{aligned}
    }_{S_2}\\ &+\underbrace{\sum_{1\le|k|\le q_{n}-1}}_{S_3}+\underbrace{\sum_{\substack{q_n+1\le|k|\le \tau_nq_{n} \\ k\ne jq_n}}}_{S_4}+\underbrace{\sum_{|k|\ge \tau_nq_n}}_{S_5}+\underbrace{\sum_{\substack{|k|=jq_n\\2\le j \le \tau_n}}}_{S_6}\end{split}
\end{equation*}
In the estimation of $S_4, S_5$ we will explicitly take $\tau_n=q_{n+1}/4q_n$.

i) Estimation of $S_1$ and $S_2$

We know that $\hat{f}_{q_n}=\overline{\hat{f}_{-q_n}}$ because $f$ is a real function, so 
\begin{equation*}
    S_1=2\pi iq_n \cdot 2i\operatorname{Im}(\hat{f}_{q_n}e^{i2\pi q_nx})|X(m,q_n)|
\end{equation*}

To estimate the term $\operatorname{Im}(\hat{f}_{q_n}e^{i2\pi q_nx})$ from below, we define a good set $I_n$ for $n\ge 4$
\begin{equation}\label{I_n}
    I_n:=\bigg\{x\in\mathbb{T}:\{ q_nx+\frac{\arg(\hat{f}_{q_n})}{2\pi}\in\bigg[\frac{1}{n},\frac{1}{2}-\frac{1}{n}\bigg]\bigcup\bigg[\frac{1}{2}+\frac{1}{n},1-\frac{1}{n}\bigg]\}\bigg\}
\end{equation}

When $x\in I_n$, $|\operatorname{Im}(\hat{f}_{q_n}e^{i2\pi q_nx})|=|\hat{f}_{q_n}||\sin{(2\pi q_n x+\arg \hat{f}_{q_n})}|\ge 4|\hat{f}_{q_n}|/n$. Moreover, $|I_n|\to1$, so it satisfies the requirement of criterion \ref{Criterion of mixing along a sequence}.
Then since $q_n\le m\le q_{n+1}/2$, using estimate \ref{6} we obtain
\begin{equation*}
    |S_1|=2\pi q_n |X(m,q_n)||\operatorname{Im}(\hat{f}_{q_n}e^{i2\pi q_nx})|\ge c_1|\hat{f}_{q_n}|\cdot m\cdot{q_n}/n
\end{equation*}
For $S_2$, we have the inequality: 
\begin{equation*}
    |X(m,q_n)-|X(m,q_n)||\le2|X(m,q_n)|\arg(X(m,q_n))
\end{equation*}
combined with the trivial bound \ref{1} and the estimate \ref{7}:
{
\begin{equation*}
\begin{split}
    |S_2|=2\pi q_n \cdot |\hat{f}_{q_n}|\bigg|X(m,q_n)-|X(m,q_n)|\bigg|&\le 2\pi q_n \cdot |\hat{f}_{q_n}|\cdot 2m\cdot \pi\frac{m}{q_{n+1}}\\&= c_2|\hat{f}_{q_n}|\cdot m^2\cdot q_n/q_{n+1}
\end{split}
\end{equation*}
}

ii) Estimation of $S_3$ and $S_4$
\enlargethispage{1\baselineskip}
{
Direct computation shows
\begin{equation*}
    |(S_3+S_4)|\le \sum_{1\le|k|\le q_{n}-1}2\cdot2\pi k\cdot |\hat{f}_k||X(m,k)|+\sum_{\substack{q_n+1\le|k|\le \tau_nq_{n} \\ k\ne jq_n}}2\cdot2\pi k\cdot |\hat{f}_k||X(m,k)|
\end{equation*}
using \ref{2} and \ref{3} we obtain
\begin{equation*}
\begin{split}
     |(S_3+S_4)|&\le \sum_{k=1}^{q_n-1}2\cdot2\pi k\cdot |\hat{f}_k|\cdot q_n + \sum_{k=q_n+1}^{\tau_nq_n}2\cdot2\pi k\cdot |\hat{f}_k|\cdot2q_n\\&=C_3q_n\sum_{k=1}^{q_n-1}k\cdot |\hat{f}_k|+C_4q_n\sum_{k=q_n+1}^{\tau_nq_n}k\cdot |\hat{f}_k|
\end{split}
\end{equation*}
Since $f$ is analytic, there exists $0<r<1$ such that $|\hat{f}_k|<Cr^{|k|}$. Therefore $\sum k\cdot |\hat{f}_k|$ is bounded. We can then find a constant $c_3$ such that \begin{equation*}
    |(S_3+S_4)|\le c_3q_n
\end{equation*}}

iii) Estimation of $S_5$

Again $\tau_n=q_{n+1}/4q_n$. Using the naive bound \ref{1} and $|\hat{f}_k|<Cr^{|k|}$, we have
\begin{equation*}
    |S_5|\le 2\cdot2\pi m\sum_{k=\tau_nq_n}^\infty k\cdot |\hat{f}_k| \le c_4m \cdot\tau_nq_nr^{\tau_nq_n}=c_5 m \cdot q_{n+1}r^{q_{n+1}/4}
\end{equation*}

iv) Estimation of $S_6$

We will again use the naive bound  \ref{1} to get
\begin{equation*}
    |S_6|\le c_6m \sum_{\substack{|k|=jq_n\\2\le j \le \tau_n}} |k\cdot \hat{f}_k|
\end{equation*}

Combining the above estimates we obtain the following:
 \begin{equation}\label{main estimate on f'}
 \begin{split}
     |S^\prime_m(f)(x)|&\ge c_1|\hat{f}_{q_n}| m\cdot{q_n}/n-c_2|\hat{f}_{q_n}|m^2\cdot q_n/q_{n+1}-c_3q_n\\&-c_5 m\cdot q_{n+1}r^{q_{n+1}/4}-c_6m \sum_{\substack{|k|=jq_n\\2\le j \le \tau_n}} |k|\cdot |\hat{f}_k|
 \end{split}
\end{equation}

There is no uniform bound for the last term without further information of $\hat{f}_k$, but the following lemma \ref{Fourier estimate} shows that we can always pick a subsequence $\hat{f}_{l_k}$ to make the last term negligible {when $f$ is not a trigonometric polynomial.}
\begin{lem}\label{Fourier estimate}
    Let $f$ be a real analytic function that is not a trigonometric polynomial. Let its  Fourier expansion be $f=\sum \hat{f}_k e^{2\pi i kx}$,  and write $T_{n}=\sum_{|j|\ge2}|\hat{f}_{j n}|\cdot |j n|$. Then there exists a subsequence $\{\hat{f}_{l_k}\}_{k\ge 1}$ of the Fourier coefficients such that 
    \begin{equation}\label{T_l_k estimate}
        \lim_{k\to\infty}\bigg|\frac{\hat{f}_{l_k}}{T_{l_k}}\bigg|\to \infty
    \end{equation}
\end{lem}

\begin{proof}
Since $f$ is real, $|\hat{f_k}|=|\hat{f_{-k}}|$, thus $T_{l_k}=2\sum_{j\ge2}|\hat{f}_{j n}|\cdot |j n|$. 
  Since $f$ is analytic, the power series $P(z)=\sum_{n=0}^\infty \hat{f}_nz^n$ has radius of convergence $R>1$. For $1<\rho<R$, $|\hat{f}_n|\rho^n\to 0$. Define $u(\rho)=\max_{n\ge 0}|\hat{f}_n|\rho^n$, and $v(\rho)=\max\{n\in\mathbb{N}: |\hat{f}_n|\rho^n=u(\rho)\}$. Then we choose a sequence of strictly increasing real numbers $\{\rho_k\}_{k\ge 1}$ such that $1<\rho_k<R$ and $\rho_k\to R$. Since $f$ is not a trigonometric polynomial, such a sequence always exist. We let $l_k=v(\rho_k)$. By definition of $v(\rho_k)$, for any $j\ge 2$, we have the inequality 
  \begin{equation}
      |\hat{f}_{l_k}|\rho_k^{l_k}\ge|\hat{f}_{jl_k}|\rho_k^{jl_k}\implies |\hat{f}_{jl_k}|\le |\hat{f}_{l_k}|\rho_k^{-l_k(j-1)}
  \end{equation}
  Substitute this into the expression for $T_{l_k}$:
  \begin{equation}
      T_{l_k}\le2\sum_{j\ge2}|\hat{f}_{l_k}|\rho_k^{-l_k(j-1)}\cdot jl_k=2|\hat{f}_{l_k}|l_k\sum_{j\ge2}\rho_k^{-l_k(j-1)}\cdot j
  \end{equation}
  using 
  \begin{equation*}
      \sum_{j\ge 2}x^{j-1}j=\frac{d}{dx}\bigg(\sum_{j\ge 0}x^{j}-x-1\bigg)=\frac{1}{(1-x)^2}-1=\frac{2x-x^2}{(1-x)^2}
  \end{equation*}
  when $k\to\infty$ and therefore $x\to 0$, RHS is bounded by $4x$, so
  \begin{equation}
      \bigg|\frac{\hat{f}_{l_k}}{T_{l_k}}\bigg|\ge\frac{|\hat{f}_{l_k}|}{2|\hat{f}_{l_k}|l_k\cdot4\rho_k^{-l_k}}=\frac{\rho_k^{l_k}}{8l_k}\xrightarrow[k\to\infty]{}\infty
  \end{equation}
  and this completes the proof.
\end{proof}

On the other hand, since $f\in C^2(\mathbb{T})$, we have the following naive estimate on the second derivative: 
\begin{lem}\label{upper bound for f''}
    There exists a constant $C$ such that $\forall m$ and $x\in\mathbb{T}$, 
    \begin{equation*}
        \bigg|S^{\prime\prime}_m(f)(x)\bigg|\le Cm
    \end{equation*}
\end{lem}

\section{Proof of theorem \ref{main thm}}

The proof of theorem \ref{main thm} involves two steps: For any irrational $\beta$ , we first construct a dense $G_\delta$ set $\mathcal{U}\in\mathbb{T}$ based on the convergents of $\beta$, denoted by $\{q_n'\}_{n\ge1}$, and then show that for every $\alpha\in\mathcal{U}$, $T_\alpha^f$ is mixing along a subsequence $\{q_{\eta(n)}'\}_{n\ge 1}$. Theorem \ref{main thm} then follows by applying Theorem \ref{main criterion} and lemma \ref{lemma on rigid sq} by letting $d_n=q_{\eta(n)}'$ .

Denote the convergents of $\beta$ by $\{q_n'\}_{n\ge1}$. For an analytic function $f$ that is not a trigonometric polynomial, let $L=\{l_k\}_{k\ge 1}$ be the set given by the lemma \ref{Fourier estimate}. For $k\ge1$, define $\eta(k)$ as the smallest integer such that \begin{equation}\label{q_eta(k)' estimate}
   q_{\eta(k)}'>\frac{2}{|\hat{f}_{l_k}|^2}
\end{equation} 
Notice that the lemma \ref{Fourier estimate} implies $\hat{f}_{l_k}^2\ne0$. 
Define 
    \begin{equation*}
        V_{n,m}:=\bigg\{\alpha:q_n(\alpha)=l_m, \text{and }q_{n+1}>2q_{\eta(m)}'\cdot n^2\bigg\}
    \end{equation*}
    and then let 
    \begin{equation*}
\mathcal{U}_j:=\bigcup_{n=1}^\infty V_{n,j}
    \end{equation*}
$\mathcal{U}_j$ is the set of $\alpha$ that satisfies:

    i) Among the convergents of $\alpha$,  there exists a $n\in\mathbb{N} $ s.t $q_n(\alpha)= l_j$,

ii) For $n$ we pick in i), we have
\begin{equation*}
    q_{n+1}>2{q_{\eta(j)}'}\cdot n^2
\end{equation*}
\begin{lem}\label{U_j open}
    Each $\mathcal{U}_j$ is open in the subspace topology on $\mathbb{R}\setminus\mathbb{Q}$. 
\end{lem}
\begin{proof}
     Consider in general a sequence of partial quotient $\omega=[c_0;c_1,c_2\cdots]$. There is an unique sequence of convergents $q_n(\omega)$ associated to $\omega$ that is determined by $\omega$ and the recursive formula \ref{recursive formula continued fraction}. (For example, $q_1(\omega)=c_1$ and $q_2(\omega)=c_1q_1(\omega)+1$, etc.).
     For $\omega\in V_{n,m}$, define \begin{equation*}
    A_\omega:=\bigg\lceil\frac{2q_{\eta(m)}'\cdot n^2-q_{n-1}(\omega)}{l_m}\bigg\rceil
     \end{equation*}Then $q_{n+1}(\omega)>2q_{\eta(m)}'\cdot n^2$ if and only if $c_{n+1}\ge A_\omega$.
     For fixed $m$, there are only finitely many tuples $\omega=(c_1,c_2\cdots c_n)$ such that $q_{n}(\omega)=l_m$.
    
     We can then rewrite $V_{n,m}$ as an union of cylinder sets
     \begin{equation*}\label{V_n,m decomposition}
        V_{n,m}=\bigcup_{c_0\in\mathbb{Z}}\bigcup_{\substack{\omega=(c_1,\cdots c_n) \\\text{s.t }q_n=l_m }}\bigcup_{c_{n+1}=A_\omega}^\infty\bigg(I(c_0;c_1,\cdots c_n,c_{n+1})\cap(\mathbb{R}\setminus\mathbb{Q})\bigg)
     \end{equation*}
     where $I(c_0;c_1,\cdots c_n,c_{n+1})$ is the set that contains all irrational number whose partial quotients start with $(c_0;c_1,\cdots c_n,c_{n+1})$, which is open in the subspace topology on $\mathbb{R}\setminus\mathbb{Q}$. Therefore, $V_{n,m}$ and $\mathcal{U}_j$ are both open in the subspace topology on $\mathbb{R}\setminus\mathbb{Q}$ for all $n,m\ge 1$. 
\end{proof}

\begin{lem}\label{U_j dense}
    Each $\bigcup_{j=M}^\infty\mathcal{U}_j$ is dense on $\mathbb{R}$ for all $M\ge 1$.
\end{lem}
\begin{proof}
We need the following Legendre's Criterion in continued fractions to determine if a given rational number $p/q$ is a convergent of an irrational number $\alpha$.
\begin{lem}(Legendre's Criterion)
Let $\alpha \in \mathbb{R}$, $p,q \in \mathbb{Z}^+$ with $\gcd(p, q) = 1$. If 
\[
\left| \alpha - \frac{p}{q} \right| < \frac{1}{2q^2}
\]
then $\frac{p}{q}$ is a convergent of the continued fraction expansion of $\alpha$.
\end{lem}
For any interval $(c,d)\in\mathbb{R}$, we take $j$ so large such that 
\begin{equation*}
    \frac{1}{l_j}< \frac{1}{10(d-c)}
\end{equation*}
then $\forall k>j$ there exists $m_k\in\mathbb{N}$ where $\gcd(m_k, l_k)=1$, and
\begin{equation*}
    c+\frac{3}{10}\cdot \frac{1}{d-c}<\frac{m_k}{l_k}<d-\frac{3}{10}\cdot \frac{1}{d-c}
\end{equation*}
Now consider the interval 
\begin{equation*}
    J_k=\bigg(\frac{m_k}{l_k}-\frac{1}{2l_k^2},\frac{m_k}{l_k}+\frac{1}{2l_k^2}\bigg)
\end{equation*}
Clearly, $J_k\subset(c,d)$. By Legendre's criterion, for any $\alpha\in J_k\cap (\mathbb{R}\setminus\mathbb{Q})$, $m_k/l_k$ is a convergent of the continued fraction expansion of $\alpha$. Let the partial quotient of $m_k/l_k$ be $\omega=[c_0;c_1,\cdots,c_n]$, we may take $c_{n+1}$ so large that any irrational number $\alpha$ whose partial quotients start with $[c_0;c_1,\cdots,c_n,c_{n+1}]$ is in $V_{n,j}\subset\mathcal{U}_j$. Moreover, if $c_{n+1}$ is large enough, then $\alpha\in J_k\subset(c,d)$, which implies that $\mathcal{U}_j\cap(c,d)\neq\varnothing$.
Therefore, $\bigcup_{j=M}^\infty\mathcal{U}_j$ is dense.
\end{proof}

Finally, we define $\mathcal{U}$ as the lim sup of $\mathcal{U}_j$:
\begin{equation}\label{intersection decomposition}
\mathcal{U}=\bigcap_{M=1}^\infty\bigcup_{j=M}^{\infty}\mathcal{U}_j
\end{equation}
In plain words, $\mathcal{U}$ is the set of irrational numbers $\alpha$ so that there is a subsequence $q_{\sigma(n)}$ where $q_{\sigma(n)}(\alpha)\in L$ for every $n$. To simplify the notation, we delete all elements of $L$ that are not of  form $q_{\sigma(n)}(\alpha)$, then write $q_{\sigma(n)}(\alpha)=l_n$. This new set still has all properties we constructed above, in particular,
\begin{equation}\label{q_[sigma(n)+1] estimate}
q_{\sigma(n)+1}>2q_{\eta(n)}'\cdot\sigma(n)^2.
\end{equation}

\begin{thm}\label{Gdelta dense thm}
    $\mathcal{U}$ is $G_\delta$ dense in $\mathbb{R}$.
\end{thm}

\begin{proof}
    
 By lemma \ref{U_j open} and \ref{U_j dense}, $\bigcup_{j=M}^\infty\mathcal{U}_j$ is both open and dense for all $M\ge1$, therefore, by the Baire Category theorem, $\mathcal{U}$ is $G_\delta$ dense. 
\end{proof}

\begin{thm}\label{mixing subsq thm}
    Let $\{q_n'\}_{n\ge1}$ be the convergents of $\beta$, then for every $\alpha\in\mathcal{U}$ the flow $T_\alpha^f$ is mixing along the subsequence $\{q_{\eta(n)}'\}_{n\ge 1}$.
\end{thm} 
\begin{proof}
We adopt the notation $f(n)\gg g(n)$ if there exists a constant $C$ such that $f(n)>Cg(n)$ for sufficiently large $n$. To simplify the notation, in the rest of the section we write $b_n := \hat{f}_{q_{\sigma(n)}}$ as the corresponding Fourier
coefficient of $f$. By the definition of $\eta(n)$, see \ref{q_eta(k)' estimate}, for large $n$ we have $q_{\eta(n)}'>2/|\hat{f}_{l_n}|^2\gg l_n=q_{\sigma(n)}$, and inequality \ref{q_[sigma(n)+1] estimate} gives $q_{\sigma(n)+1}/2>2q_{\eta(n)}'$. Therefore, for $m_n\in [q_{\eta(n)}'/2,2q_{\eta(n)}']$, we can apply the estimate \ref{main estimate on f'} and Lemma \ref{Fourier estimate} to obtain the following bound on $|S_{m_n}^\prime(f)|$:
\begin{equation}
     \begin{split}
     |S^\prime_{m_n}(f)(x)|&\ge \underbrace{c_1|b_{n}| m_n\cdot{q_{\sigma(n)}}/\sigma(n)}_{S_1}-\underbrace{c_2|b_{n}|m_n^2\cdot q_{\sigma(n)}/q_{\sigma(n)+1}}_{S_2}-\underbrace{c_3q_{\sigma(n)}}_{S_3+S_4}\\&-\underbrace{c_5 m_n\cdot q_{\sigma(n)+1}r^{q_{\sigma(n)+1}/4}}_{S_5}-\underbrace{c_6m_n\cdot o(|b_{n}|)}_{S_6}
 \end{split}
\end{equation}

By inequality \ref{q_[sigma(n)+1] estimate}, 
\begin{equation*}
    \frac{S_1}{S_2}=\frac{c_1q_{\sigma(n)+1}}{c_2m_n\sigma(n)}>\frac{c_1\cdot2q_{\eta(n)}'\cdot\sigma(n)^2}{c_2\cdot 2q_{\eta(n)}'\cdot\sigma(n)}=\frac{c_1\sigma(n)}{c_2}\xrightarrow[n\to\infty]{}\infty
\end{equation*}and the inequality \ref{q_eta(k)' estimate} implies \begin{equation*}
    m_n>\frac{q_{\eta(n)}'}{2}>\frac{1}{|\hat{f}_{l_n}|^2}=\frac{1}{|b_n|^2}
\end{equation*}
therefore, using the fact that $b_n\sigma(n)\to0$, 
\begin{equation*}
    \frac{S_1}{S_3+S_4}=\frac{c_1|b_n|m_n}{c_3\sigma(n)}>\frac{c_1}{c_3\cdot\sigma(n)\cdot|b_n|}\xrightarrow[n\to\infty]{}\infty
\end{equation*}
For $S_5$, observe that $S_5\to 0$ and $S_1$ is not, so $S_1\gg S_5$. Finally, $q_{\sigma(n)}\gg \sigma(n)$ gives $S_1\gg S_6$.  Combining these estimates, we obtain:
\begin{equation}\label{Final estimate of f'}
    |S^\prime_{m_n}(f)(x)|\ge C'\cdot S_1= C^\prime|b_{n}| m_n\cdot{q_{\sigma(n)}}/\sigma(n)
\end{equation}
for some constant $C^\prime$ for large enough $n$.

We now check the criterion \ref{Criterion of mixing along a sequence}. Take $t_n=q_{\eta(n)}'$, by \ref{N(x,t) estimation}, $N(x,t_n)\in [t_n/2,2t_n]$. For the sequence $q_{\eta(n)}'$, take a partial partition of the set $I_{\sigma(n)}$ defined in the previous section with intervals $\mathcal{P}_i^n$ of length between $|b_{n}|$ and $2|b_{n}|$. Since $|I_{\sigma(n)}|>1-4/\sigma(n)\to 1$ and $|b_{n}|\to 0$, partition $\mathcal{P}_i^n$ satisfies the condition \ref{eq:M}. Moreover, 
\begin{equation*}
\begin{split}
    \min_{I \in \mathcal{P}_n} \min_{m_n \in [\underline{N}(I, t_n), \overline{N}(I, t_n)]} \min_{x \in I} |S'_{m_n}(f)(x)| |I| &\geq C^\prime |b_{n}|^2 \cdot \underline{N}(I,t_n)\cdot{q_{\sigma(n)}}/\sigma(n)\\&>\frac{C'q_{\sigma(n)}}{\sigma(n)}\xrightarrow[n\to\infty]{} \infty
\end{split}
\end{equation*}
the last inequality uses the inequality \ref{q_eta(k)' estimate}. Then applying \ref{upper bound for f''} and \ref{Final estimate of f'} we have 
\begin{equation} 
\begin{split}
    \max_{I \in \mathcal{P}_n} \frac{\max_{0 \le m_n \le 2t_n} \max_{x \in I} |S''_{m_n}(f)(x)| |I|}{\min_{m_n \in [\underline{N}(I, t_n), \overline{N}(I, t_n)]} \min_{x \in I} |S'_{m_n}(f)(x)|} &\le \frac{C\cdot 2t_n\cdot  2|b_{n}|}{C^\prime|b_{n}|\cdot \underline{N}(I,t_n)\cdot{q_{\sigma(n)}}/\sigma(n)}\\&=\frac{8C}{C^\prime}\cdot\frac{\sigma(n)}{q_{\sigma(n)}}\to 0
\end{split}
\end{equation}
Thus, our construction meets all conditions of the criterion \ref{Criterion of mixing along a sequence}, which means that $T_{\alpha}^f$ is mixing along $q_{\eta({n})}'$. 
\end{proof}

\begin{proof}[proof of Theorem \ref{main thm}] 
    The main theorem now follows by applying theorem \ref{main criterion} to the set $\mathcal{U}$.
\end{proof}

\section{Von Neumann Flows with one discontinuity}\label{Von flow Chapter}
Let $\mathcal{T} = (T^f_
t )_{t\in\mathbb{R}}$ and $\mathcal{R} = (R^f_
t )_{t\in\mathbb{R}}$ be two special flows over irrational rotations $T, R : \mathbb{T} \to \mathbb{T}, T (x) = x + \alpha, R(x) = x + \beta$, under the same roof function $f:\mathbb{T} \to \mathbb{R}^+, \int_{\mathbb{T}} f= 1$, of the form
\begin{equation*}
    f (x) = A_f \{x\} +f_{ac}(x), \quad f_{ac}\in C^1{(\mathbb{T})}
\end{equation*}
The following lemma from \cite{KANIGOWSKI_SOLOMKO_2018} will allow us to neglect the $C^1$
part of the roof function. See the last paragraph in the proof of Theorem \ref{thm section 7}.

\begin{lem}\label{f_ac estimate}
    Let $g \in C^1(\mathbb{T})$. For every $\epsilon> 0$ there exists $\delta > 0$ such

that for every $x, y \in \mathbb{T}$ with $\|x-y\|<\delta$ and every $n \in\mathbb{Z}$
\begin{equation*}
    |S_n(g)(x)- S_n(g)(y)| < \epsilon \max\{1, |n|\|x-y\|\}.
\end{equation*}
\end{lem}

We will briefly explain the notation before stating the main technical lemma. Again, we denote the denominators of the convergents of $\alpha$ and $\beta$ by $\{q_n\}_{n\ge1}$ and $\{q_n'\}_{n\ge1}$. For a fixed $k$, we denote by $q'_{\sigma(k)}$ the unique number such that $q'_{\sigma(k)}\le q_k<q'_{\sigma(k)+1}$, and by $q_{\eta(k)}$ the unique number such that $q_{\eta(k)}\le q'_k<q_{\eta(k)+1}$. Then for $k=1,2,\cdots$, we define
\begin{equation}
\begin{split}
     &M_k=\max_{r:q'_{\sigma(k)}<q_r<q'_{\sigma(k)+1}}\min\bigg\{\frac{q_r}{q'_{\sigma(k)}},\frac{q'_{\sigma(k)+1}}{q_r}\bigg\},  \\&N_k=\max_{r:q_{\eta(k)}<q_r'<q_{\eta(k)+1}}\min\bigg\{\frac{q_r'}{q_{\eta(k)}},\frac{q_{\eta(k)+1}}{q_r'}\bigg\}
\end{split}
\end{equation}

We need the following theorem from the classical work of H.Kesten that describes the distribution of successive values of $\{n\alpha\}_{n\le N}$ for any fixed $N$. Let $\alpha=[a_1,a_2,\cdots]$ and its $n^{th}$ convergents $p_n/q_n$, we introduce
\begin{equation*}
    a_{n+1}^*=a_{n+1}+[a_{n+2},a_{n+3}\cdots]=a_{n+1}+\frac{1}{a_{n+2}'}
\end{equation*}
and 
\begin{equation*}
    q_{n+1}^*=a_{n+1}'q_{n}+q_{n-1}=q_{n+1}+\frac{q_{n}}{a_{n+2}'}
\end{equation*}
\begin{thm} (Theorem 1,\cite{Kesten1966})

  For any $N\in\mathbb{N}$, the Ostrowski expansion of $N$ is given by $N=\sum c_iq_i=\sum_{i=0}^{m(N,\alpha)}c_i(N,\alpha)q_i(\alpha)$ where $0\le c_i\le a_{i+1}$. For each $0\le r\le q_m-1$, the interval $(r/q_m,(r+1)/q_m)$  contains exactly one point $\{k\alpha\}$ with $1\le k\le q_m$. Denote the point in $(r/q_m,(r+1)/q_m)$ by $P_r$ and the interval $[P_r,P_{r+1})$ by $J_r$. Then exactly $q_m-q_{m-1}$ intervals $J_r$ have length $a_{m+1}^*/q_{m+1}^*$ and exactly $q_{m-1}$ have length $(a_{m+1}^*+1)/q_{m+1}^*$. The next $(c_m-1)q_m$ points subdivide the intervals $J_r$, in such a manner that
exactly $(c_m-1)$ points fall in each $J_r$. Of the last $N -c_mq_m$ points $\{k\alpha\}$, $c_mq_m+1\le k\le N$, at most one will belong to each $J_r$. These points divide each $J_r$ into $\Phi_r= c_m, \textit{or }c_m+1$ sub-intervals. Starting from $P_r$, the first $\Phi_r-1$ sub-intervals of $J_r$ have length $1/q_{m+1}^*$ and the last sub-interval has length $|J_r|-(\Phi_r-1)/q_{m+1}^*$.  
\end{thm}

Now we state the main theorem for this section:
\begin{thm}\label{thm section 7}
    If there exists a sequence $\{n_j\}$ such that $\lim_{j\to\infty}M_{n_j}=\infty$ or $\lim_{j\to\infty}N_{n_j}=\infty$, then flows $\mathcal{T} = (T^f_
t )_{t\in\mathbb{R}}$ and $\mathcal{R} = (R^f_
t )_{t\in\mathbb{R}}$ are spectrally disjoint.
\end{thm}
\begin{proof}
    WLOG we may assume that $\lim_{j\to\infty}M_{n_j}=\infty$, and by slightly abusing the notation we let $q_n$ be the convergent that attains the maximum in the definition of $M_n$.  By the Denjoy-Koksma inequality, for every $x\in\mathbb{T}$
\begin{equation*}
    \bigg|S_{q_n}(f)(x)-q_n\int f\bigg|<Var(f)=A_f
\end{equation*}
therefore, for the flow $\mathcal{T} = (T^f_
t )_{t\in\mathbb{R}}$, if we take $h_n=q_n$ and $d_n=q_n\int f=q_n$ then $S_{h_n}f(x)-d_n$ has exponential tails. By our main criterion \ref{main criterion}, we only need to show that the other flow $\mathcal{R} = (R^f_
t )_{t\in\mathbb{R}}$ is mixing along a subsequence of $q_n$.

Now we construct a partial partition $\mathcal{P}_n$ by cutting $\mathbb{T}$ into intervals with endpoints in $\{k\alpha\}_{1\le k\le 2n}$ and then skipping all shortest intervals. More precisely, we first cut $\mathbb{T}$ into $q_{\sigma(n)}$ intervals of $\{J_r\}_{0\le r\le q_{\sigma(n)}-1}$ and then remove the first $\Phi_r-1=\lfloor{2q_n/q_{\sigma(n)}'}\rfloor \textit{ or }\lfloor{2q_n/q_{\sigma(n)}'}\rfloor+1$ sub-intervals within each $J_r$. We know that $\max_{I\in \mathcal{P}_n}|I|\to 0 $ and the total size of $\mathcal{P}_n$ is
\begin{equation*}
    Leb (\mathcal{P}_n)>1- (2q_n-q_{\sigma(n)}')\cdot\frac{1}{q_{\sigma(n)+1}'}>1-\frac{2q_n}{q_{\sigma(n)+1}}>1-\frac{2}{M_n}\to 1, 
\end{equation*}
By proposition \ref{criterion of mixing by uniform stretch} and lemma \ref{N(x,t) estimation}, it is now sufficient to show that for any $r \in [q_n/2,2q_n]$, any $x\in\mathbb{T}$ and any $I\in\mathcal{P}_n$, $S_r(f)(x)$ is $(\epsilon(n),k(n))$ uniformly stretching on $I$ for some $\epsilon(n)\to 0$ and $k(n)\to\infty$.

 Ignore the $C^1$ part for now, then by our construction, for any $I\in \mathcal{P}_{n}$, $T^j(I)\cap \{0\}=\varnothing, \forall 0<j<r<2q_n$. In other words, $T^r(I)$ will not hit the discontinuity. Each interval $I\in\mathcal{P}_n$ satisfies:
\begin{equation*}
\begin{split}
    |I|=|J_r|-(\Phi_r-1)\cdot\frac{1}{q_{\sigma(n)+1}'}&\ge\frac{a_{\sigma(n)+1}'-(\lfloor{2q_n/q_{\sigma(n)}'}\rfloor+1)}{q_{\sigma(n)+1}'}\\&>\frac{a_{\sigma(n)+1}'(1-2/M_n)}{q_{\sigma(n)+1}'}
\end{split}
\end{equation*}
Recall that $q_{\sigma(n)+1}'<(a_{\sigma(n)+1}'+1)q_{\sigma(n)}$.  Therefore, for large $n$ the total stretch on each $I$ is bounded by:
\begin{equation*}
\begin{split}
    \sup_{x\in I}S_r(f)(x)-\inf_{x\in I}S_r(f)(x)=A_f\cdot r|I|&> A_f\cdot\frac{q_n}{2}\frac{a_{\sigma(n)+1}'(1-2/M_n)}{q_{\sigma(n)+1}'}\\&>A_f\cdot M_n\cdot\frac{1}{8}
\end{split}
\end{equation*}
which tends to infinity by our assumptions.

To show that the stretch is uniform, observe that in each interval $I$ the contribution from $A_f\{x\}$ part is linear, and by lemma \ref{f_ac estimate} the contribution from $C^1$ part oscillates by $\epsilon r|x-y|$ for $x,y$ in any sufficiently small intervals, depending on $\epsilon$. Therefore, for any function $\epsilon(j)\to0$, we can pick a sufficiently large $n_j$, where intervals in $\mathcal{P}_{n_j}$ are less than $\delta_{\epsilon(j)}$ as in \ref{f_ac estimate}. In this way, we guarantee that for $I\in\mathcal{P}_{n_j}$, $S_r(f)$ is $(2\epsilon(j), A_fM_{n_j}/8)$ uniformly stretching, which in turn implies that $\mathcal{R} = (R^f_
t )_{t\in\mathbb{R}}$ is mixing along $q_{n_j}$ as a subsequence of $q_n$.
\end{proof}

\begin{section}{The set of $\beta$ is of full measure}\label{sec beta full measure}

\begin{thm}\label{thm full measure set}
    For every $\alpha\in[0,1]\setminus\mathbb{Q}$, there exists a full measure set $S_\alpha$, $m(S_\alpha)=1$, such that for every $\beta\in S_\alpha$,  there exists a sequence $\{n_j\}_{j\ge1}$ such that either $\lim_{j\to\infty}M_{n_j}=\infty$ or $\lim_{j\to\infty}N_{n_j}=\infty$.
\end{thm}
We need the following technical lemma:
\begin{lem}\label{a_n' lemma}
    Let $\alpha=[a_1,a_2,\cdots]$ and $\beta=[a_1',a_2',\cdots]$. If both $M_k$ and $N_k$ are bounded by a constant $C$, then there exists an absolute constant $K$, such that for all $n\in\mathbb{N}$, there exists $1\le j\le \lfloor K\ln C\rfloor$ s.t
    \begin{equation}\label{a_n' upper and lower bound}
        \frac{1}{C}\cdot a_{\eta(n)}a_{\eta(n)+1}\cdots a_{\eta{(n)}+j}-1<a_n'<C(a_{\eta(n)}+1)(a_{\eta(n)+1}+1)\cdots (a_{\eta{(n)+j}}+1)
    \end{equation}
\end{lem}
\begin{proof}
    We consider the following sequence of numbers:
    \[q_{\eta(n)}<q_n'<q_{\eta(n)+1}<\cdots<q_{\eta(n+1)}<q_{n+1}'<q_{\eta(n+1)+1}\]
    It is possible that $q_{\eta(n+1)}=q_{\eta(n)+1}$. Since $M_n<C$ and convergents grow at least exponentially, there are at most $\lfloor K\ln C\rfloor$ terms between $q_n'$ and $q_{n+1}'$, where $K$ is an absolute constant independent of any other quantity. Since $N_n<C$, by definition we have 
    \begin{equation*}
        \min\bigg\{\frac{q_n'}{q_{\eta(n)}},\frac{q_{\eta(n)+1}}{q_n'}\bigg\}<C, \text{\quad and }\min\bigg\{\frac{q_{n+1}'}{q_{\eta(n+1)}},\frac{q_{\eta(n+1)+1}}{q_{n+1}'}\bigg\}<C 
    \end{equation*}

There are four possible cases but the proof are very similar. We will prove for the case that ${q_n'}/{q_{\eta(n)}}<C$ and ${q_{n+1}'}/{q_{\eta(n+1)}}<C$, and the other cases will follow from same argument. By the recursive formula of convergents, we have 
\begin{equation*}
    a_n'<\frac{q_{n+1}'}{q_{n}'}<\frac{Cq_{\eta(n+1)}}{q_{\eta(n)}}<C(a_{\eta(n)}+1)(a_{\eta(n)+1}+1)\cdots (a_{\eta{(n+1)-1}}+1)
\end{equation*}
and 
\begin{equation*}
    a_n'+1>\frac{q_{n+1}'}{q_{n}'}>\frac{q_{\eta(n+1)}}{Cq_{\eta(n)}}>\frac{1}{C}\cdot a_{\eta(n)}a_{\eta(n)+1}\cdots a_{\eta{(n+1)}-1}
\end{equation*}
which is exactly inequality \ref{a_n' upper and lower bound}.
\end{proof}
\begin{proof}[proof of Theorem \ref{thm full measure set}]
    Let $\alpha=[a_1,a_2,\cdots]$, $\beta=[a_1',a_2',\cdots]$, and assume that both $M_k$ and $N_k$ are bounded by a constant $C$. We will show that such $\beta$ form a Lebesgue measure zero set. 
    
    Fix an $n\in\mathbb{N}$, define interval $L_n^j\subset\mathbb{R}$ as:
    \begin{equation}
        L_n^j:=\bigg[\frac{1}{C}\cdot a_{\eta(n)}a_{\eta(n)+1}\cdots a_{\eta{(n)}+j}-1,C(a_{\eta(n)}+1)(a_{\eta(n)+1}+1)\cdots (a_{\eta{(n)+j}}+1)\bigg]
    \end{equation}  by the lemma \ref{a_n' lemma}, there exists $1\le j< \lfloor K\ln C\rfloor$ s.t $a_n'\in L_n^j$. For each $n$, let $A_n$ be the set of allowed value for $a_n'$, then $A_n\subset \bigcup_{j=1}^{\lfloor K\ln C\rfloor} L_n^j$. To simplify the notation, we now write $L(n,j)$ as the lower end of $L_n^j$ and $R(n,j)$ as the upper end of $L_n^j$.
    Then the length of $L_n^j$ can be estimated by 
    \begin{equation}
    \begin{split}
        |L_n^j|&<C^2\cdot (L(n,j)+1)\cdot\frac{(a_{\eta(n)}+1)(a_{\eta(n)+1}+1)\cdots (a_{\eta{(n)+j}}+1)}{a_{\eta(n)}a_{\eta(n)+1}\cdots a_{\eta{(n)}+j}}\\&<C^22^{\lfloor K\ln C\rfloor}L(n,j)
    \end{split}
    \end{equation}

Let $x=[k_0; k_1,k_2,\cdots ]$ be the continued fraction expansion of a random number $x$ that is uniformly distributed in $(0,1)$, the Gauss–Kuzmin theorem states that 
\begin{equation*}
    \lim_{n\to\infty} \mathbb{P}(k_n=k)=-\log_2\bigg(1-\frac{1}{(k+1)^2}\bigg)
\end{equation*}
As a consequence, we now show that the measure of $\beta\in (0,1)\setminus\mathbb{Q}$ such that $a_n'\in A_n$ is bounded away from $1$. 
\begin{lem}
    There exists an uniform constant $\rho$ such that 
    \begin{equation*}
        \lim_{n\to\infty}\bigg| \mathbb{P}(k_n\in A_n)\bigg|<1-\rho
    \end{equation*}
\end{lem}
\begin{proof}
The idea is to show that $L_n^j$ can not cover all positive integers. We know that $L(n,1)<L(n,2)<\cdots L(n,\lfloor K\ln C\rfloor)$. If $L_n^j$ does not miss out any integers, then $R(n,j)$ is bounded by the following sequence:
\begin{equation*}
    \Sigma_0=1, \quad \Sigma_j=C^22^{\lfloor K\ln C\rfloor}(\Sigma_{j-1}+1)
\end{equation*}
Since $\lfloor K\ln C\rfloor$ is bounded, $|\bigcup_{j=1}^{\lfloor K\ln C\rfloor} L_n^j|$ is bounded. On the other hand, if there exists $1\le j< K\ln C$ such that $R(n,j)< L(n,j+1)-1$, then there exists $m<\Sigma_{\lfloor K\ln C\rfloor}$ such that $m\notin \bigcup_{j=1}^{\lfloor K\ln C\rfloor} L_n^j$. 
This shows that $\cup L_n^j$ always miss out a set of positive Gauss measure. More precisely, 
\begin{equation*}
    \lim_{n\to\infty}\bigg| \mathbb{P}(k_n\in A_n)\bigg|<1-\min\bigg\{\mathbb{P}(k_n=m),\sum_{s>S_{C^22^{K\ln C}}}\mathbb{P}(k_n=s)\bigg\}
\end{equation*}
let $\rho$ be the minium of the two terms then we are done. 
\end{proof}

    The lemma implies that there exists a positive integer $N$ such that $\forall n\ge N$, the allowed value of $a_n'$ misses a set of positive Gauss measure. Consequently, the set of $\beta$ is of measure zero.
\end{proof}
\begin{proof}[proof of Theorem \ref{Von main thm}]
     The theorem follows by combining theorem \ref{thm section 7} and theorem \ref{thm full measure set}.
\end{proof}

\end{section}
\section{acknowledgement}
I would like to thank Professor Adam Kanigowski for introducing me to this problem, and for his insightful feedback and valuable comments on this paper.

\bibliographystyle{plain} 
\bibliography{refs}
\end{document}